\documentclass[12pt]{amsart}
\newtheorem{theorem}{Theorem}
\newtheorem{proposition}[theorem]{Proposition}
\newtheorem{lemma}[theorem]{Lemma}
\newtheorem{corollary}[theorem]{Corollary}

\def\eps{\varepsilon}
\def\CC{\Bbb C}

\def\ds{\displaystyle}

\def\dist{\operatorname{dist}}
\def\ord{\operatorname{ord}}

\title
{Estimates for invariant metrics on $\Bbb C$-convex domains}
\author{Nikolai Nikolov, Peter Pflug, W\l odzimierz Zwonek}

\address
{Institute of Mathematics and Informatics\\ Bulgarian Academy of
Sciences\\1113 Sofia, Bulgaria}\email{nik@math.bas.bg}

\address{Carl von Ossietzky Universit\"at Oldenburg\\
Fachbereich Mathe\-ma\-tik\\ Postfach 2503\\ D-26111 Oldenburg,
Germany}\email{pflug@mathematik.uni-oldenburg.de}

\address{Instytut Matematyki, Uniwersytet Jagiello\'nski, \L ojasiewicza 6,
30-348 Krak\'ow, Poland}\email{Wlodzimierz.Zwonek@im.uj.edu.pl}
\begin{document}



\begin{thanks}{This paper was written during the stay of the first-named
author at the Carl von Ossietzky Universit\"at Oldenburg,
(November-December 2008) supported by the DFG grant
436POL113/103/0-2. The third-named was supported by the research
grant No. N N201 361436  of the Polish Ministry of Science and
Higher Education.}
\end{thanks}

\begin{abstract} Geometric lower and upper estimates are obtained
for invariant metrics on $\Bbb C$-convex domains containing no
complex lines.
\end{abstract}

\maketitle

\section{Introduction and results}

Let $\Bbb D\subset\Bbb C$ be the unit disc. For a domain
$D\subset\Bbb C^n$ the Carath\'eodory and Kobayashi
(pseudo)metrics are defined in the following way (cf. \cite{JP}):
$$\gamma_D(z;X)=\sup\{|f'(z)X|:f\in\mathcal O(D,\Bbb
D),f(z)=0\},$$
$$\kappa_D(z;X)=\inf\{\alpha\ge 0:\exists\varphi\in\mathcal O(\Bbb
D,D):\varphi(0)=z,\alpha\varphi'(0)=X\}.$$ It is clear that
$\gamma_D\le\kappa_D.$

Recall that a domain $D\subset\Bbb C^n$ is called {\it $\Bbb
C$-convex} if any non-empty intersection with a complex line is a
simply connected domain (cf. \cite{APS,Hor}). A consequence of the
fundamental Lempert theorem is the equality $\gamma_D=\kappa_D$
for any convex domain and any bounded $\Bbb C$-convex domain $D$
with $C^2$ boundary; for the last statement use that such a domain
can be exhausted by smooth bounded strictly $\CC$-convex domains
(see \cite{Jac}).

A domain $D\subset\Bbb C^n$ is said to be {\it linearly convex}
(respectively, {\it weakly linearly convex}) if for any $a\in\Bbb
C^n\setminus D$ (for any $a\in\partial D$) there exists a complex
hyperplane through $a$ which does not intersect $D.$

Recall that the following implications hold:

\centerline{$\Bbb C$-convexity $\Rightarrow$ linear convexity
$\Rightarrow$ weak linear convexity.}

Moreover, these three notions coincide in the case of $C^1$-smooth
domains in dimension greater than 1 (cf. \cite{APS,Hor}).

{\bf (A)} For $\CC$-convex domains we shall prove the following
results for the boundary behavior of the Carath\'eodory and
Kobayashi metrics.

\begin{proposition}\label{1} Let $D$ be a $\Bbb C$-convex domain
containing no complex line through $z\in D$ in direction of $X.$
Then
$$\frac{1}{4}\le\gamma_D(z;X)d_D(z,X)\le\kappa_D(z;X)d_D(z,X)\le 1,$$
where $$d_D(z,X)=\sup\{r>0:z+\lambda X\in D\hbox{ if
}|\lambda|<r\}$$ is the distance from $z$ to $\partial D$ in
direction $X.$
\end{proposition}

The constant $\frac{1}{4}$ can be replaced by $\frac{1}{2}$ in the
case of convex domains (see \cite{BP}). On the other hand, the
constant $\frac{1}{4}$ is the best one in the planar case as the
image $D=\CC\setminus[1/4,\infty)$ of $\Bbb D$ under the Koebe
function $\frac{z}{(1+z)^2}$ shows. It is clear that the upper
constant $1$ is attained if, for example, $D=\Bbb D.$

\begin{corollary}\label{2} For any $\CC$-convex domain $D\subset\CC^n$ one has
that $\kappa_D\le 4\gamma_D.$
\end{corollary}

Recall that if a $\Bbb C$-convex domain $D\subset\Bbb C^n$
contains a complex line, then it is linearly equivalent to the
Cartesian product of $\Bbb C$ and a $\Bbb C$-convex domain in
$\Bbb C^{n-1}.$

For a boundary point $a$ of a domain $D\subset\CC^n$ denote by
$L_a$ the set of all vectors $X\in\CC^n$ for which there exists
$\varepsilon>0$ such that $\partial
D\supset\Delta_X(a,\varepsilon)=\{a+\lambda
X:|\lambda|<\varepsilon\}.$

The following result is a consequence of Proposition \ref{1}.

\begin{proposition}\label{3} Let $a$ be a boundary point of
a $\CC$-convex domain $D\subset\CC^n.$

(i) Then
$$\lim_{z\to a}\gamma_D(z;X)=\infty$$
locally uniformly in $X\not\in L_a.$\footnote{This means that for
any $M>0$ there are neighborhoods $U$ of $a$ and $V$ of $X$ such
that $\gamma_D(z;Y)>M$ for any $z\in D\cap U$ and $Y\in V\setminus
L_a.$}

(ii) If $\partial D$ is $C^1$-smooth at $a,$ then $L_a$ is a
linear space. Moreover, for any non-tangential cone $\Lambda$ with
vertex at $a$ there is a constant $c>0$ such that
$$\limsup_{\Lambda\ni z\to a}\kappa_D(z;X)\le c$$
locally uniformly in the unit vectors $X\in L_a.$
\end{proposition}

{\bf (B)} Next we shall discuss types related to a
($C^\infty$-)smooth boundary point $a$ of a domain $D\subset\CC^n$
and a vector $X\in(\CC^n)_\ast.$ Denote by $m_a$ the (usual) type
of $a,$ i.e. the maximal order of contacts of non-trivial analytic
discs through $a$ and $\partial D$ at the point $a.$ Replacing
analytic discs by complex lines, we define the linear type $l_a$
of $a.$ We may also define $l_{a,X}$ as the order of contact of
the line through $a$ in direction of $X$ and $\partial D$ at $a.$
Then $m_a\ge l_a=\sup_X l_{a,X}.$ Note that if $l_{a,X}<\infty,$
then $X\not\in L_a.$

\begin{proposition}\label{4} Let $a$ be a smooth boundary point of
a $\Bbb C$-convex domain $D\subset\CC^n$ and let $X\in
(\CC^n)_\ast$ with $l_{a,X}<\infty.$ Denote by $n_a$ the inner
normal to $\partial D$ at $a.$ Then there exist a neighborhood $U$
of $a$ and a constant $c>1$ such that
$$c^{-1}d_D(z)\le d_D(z,X)^{l_{a,X}}\le c d_D(z),\quad z\in D\cap U\cap
n_a,$$ where $d_D$ is the distance to $\partial D.$
\end{proposition}

Combining Proposition \ref{1} and \ref{4} we immediately get an
extension of the main result in \cite{Lee} from the convex to the
$\Bbb C$-convex case.

\begin{corollary}\label{5} Under the notations of Proposition
\ref{4}, there is a constant $c>0$ such that
$$c^{-1}(d_D(z))^{-1/l_{a,X}}\le\gamma_D(z;X)\le\kappa_D(z;X)
\le c(d_D(z))^{-1/l_{a,X}}.$$
\end{corollary}

The main result in \cite{McN1} (see also \cite{BS}) states that
$m_a=l_a$ for convex domains. The same remains true for a
$\CC$-convex domain.

\begin{proposition}\label{6} If $a$ is a smooth boundary point of
a $\CC$-convex domain $D\subset\CC^n,$ then $m_a=l_a.$
\end{proposition}

\noindent{\bf Remark.} We like to mention that the proof in
\cite{BS} immediately implies the above proposition in dimension
2. But we do not know if the criterion in \cite{BS} (for the
equality $m_a=l_a$) holds for any $\Bbb C$-convex domain.
\smallskip

Moreover, in the case of infinite type we have the following
result.

\begin{proposition}\label{7} If $a$ is a $C^1$-smooth boundary point of
a $\CC$-convex domain $D\subset\CC^n,$ then $\partial D$ contains
no nontrivial analytic disc through $a$ if only if $L_a=\{0\}.$
\end{proposition}

\noindent{\bf Remark.} Some of the above propositions in (A) and
(B) have local versions. In this connection recall that there is a
localization principle for the Kobayashi metric of any hyperbolic
domain (cf. \cite{JP}).
\smallskip

{\bf (C)} Now we are going to discuss multitypes of boundary
points. Recall that a smooth finite type pseudoconvex boundary
point $a$ of a domain $D\subset\Bbb C^n$ is said to be semiregular
\cite{DH} (or h-extendible \cite{Yu2}) if its Catlin multitype
$\mathcal M(a)$ coincides with its D'Angelo type $\Delta(a).$
Based on the fact that the usual type is equal to the line type in
the case of convex domains, it it is shown in \cite{Yu1} that if
$a$ is a smooth convex point (not necessarily of finite type),
then $\mathcal L(a)=\mathcal M(a)=\Delta(a),$ where $\mathcal
L(a)$ denotes the linear multitype of $a.$

We shall say that $a$ is a {\it $\Bbb C$-convex boundary point} of
a domain $D\subset\CC^n$ if there is a neighborhood $U$ of $D$
such that $D\cap U$ is a $\Bbb C$-convex domain.

\begin{proposition}\label{12} \footnote{The same result may be found in
\cite{Con}; the proof there is related on good local coordinates
and on the proof in \cite{Yu1}, whereas our proof is based on the
simple geometric Lemma \ref{11} and on the proof in \cite{Yu1}.}
If $a$ is a smooth $\Bbb C$-convex boundary point of a domain
$D\subset\CC^n,$ then $\mathcal L(a)=\mathcal M(a)=\Delta(a).$
\end{proposition}

Then the main result in \cite{Yu2} implies the following.

\begin{corollary}\label{13} Any smooth finite type $\Bbb C$-convex boundary
point $a$ of a domain  $D\subset\CC^n$ is a local (holomorphic)
peak point. Moreover, there is a neighborhood $U$ of $a$ and a
domain $\CC^n\supset G\supset\overline{D\cap U}\setminus\{a\}$
such that $a\in\partial G$ is a peak point w.r.t. the algebra
$A(G).$
\end{corollary}

This corollary is also a direct consequence of the main result in
\cite{DF}, where local holomorphic support functions which depend
smoothly on the boundary points are constructed .

We point out that the assumption of smoothness is essential as the
domain $D=\Bbb D\setminus[0,1)$ may show. It is easy to see that
the points from the deleted interval are not peak points for
$A(D).$

On the other hand, in \cite{Sib}, the following result is claimed.

\begin{proposition}\label{18} Let $D\subset\CC^n$ be a bounded convex domain.
Then $a\in\partial D$ is a peak point w.r.t. $A(D)$ if and only if
$L_a=\{0\}.$
\end{proposition}

For the convenience of the reader, we shall prove this result.

Note that there is a smooth convex bounded domain $D\subset\Bbb
C^2$ containing no non-trivial analytic discs in the boundary but
some of the boundary points (not of finite type) are not peak
points w.r.t. $A^\alpha(D)$ for any $\alpha>0$ (see \cite{Noe}).

Note also that main result in \cite{Nik1} (see also \cite{Yu3} and
\cite{BSY}) and Proposition \ref{12} give the following fact about
the boundary behavior of invariant metrics (see also
\cite{Blu,Lie}).

\begin{corollary}\label{14} Let $a$ be a finite type $\Bbb C$-convex
boundary point of a smooth bounded pseudoconvex domain $D\subset
\CC^n.$ Let $\mathcal M(a)=(m_1,\dots,m_n)$ be the Catlin
multitype of $a$ ($m_1=1$ and $m_2\le\dots\le m_n$ are even
numbers). Denote by $n_a$ the inner normal to $\partial D$ at $a.$
There is a basis $\{e_1,\dots,e_n\}$ ($e_1$ is the complex normal
vector and $\{e_2,\dots,e_n\}\subset T^\CC_a(\partial D)$) and a
constant $c>1$ such that for any $X=\sum_{j=1}^n X_je_j$ we have
$$c^{-1}\le \liminf_{n_a\ni
z\to
a}F_D(z;X)\left(\sum_j^n\frac{|X_j|}{(d_D(z))^{1/m_j}}\right)^{-1}$$
$$\le \limsup_{n_a\ni z\to
a}F_D(z;X)\left(\sum_j^n\frac{|X_j|}{(d_D(z))^{1/m_j}}\right)^{-1}\le
c.$$ Here $F_D$ is any of the Carath\'eodory, Kobayashi or Bergman
metrics.
\end{corollary}

We point out that this corollary implies Proposition \ref{4} in
the finite type case, showing in addition that for any
$X\in(\CC^n)_\ast$ there is $j=1,\dots,n$ with $l_{a,X}=m_j.$

{\bf (D)} Finally, we turn to the main part in this paper, namely,
the boundary behavior of the Bergman metric of $\Bbb C$-convex
domains. Denote by $L_h^2(D)$ the Hilbert space of all holomorphic
functions $f$ on a domain $D\subset\Bbb C^n$ that are
square-integrable and by $||f||_D$ the $L_2$-norm of $f$. Let
$K_D$ be the restriction to the diagonal to the Bergman kernel
function of $D$. It is well-known that (cf. \cite{JP})
$$K_D(a)=\sup\{|f(a)|^2:f\in L_h^2(D),\;\|f\|_D\le1\}.$$ If
$K_D(z)>0$ for some point $z\in D$, then the Bergman metric
$B_D(z;X),\\
X\in\Bbb C^n$, is well-defined and can be given by the equality
$$B_D(z;X)=\frac{M_D(z;X)}{\sqrt{K_D(z)}},$$ where
$M_D(z;X)=\sup\{|f'(z)X|:f\in L_h^2(D),\,\|f\|_D=1,\;f(z)=0\}$.

Recall that (cf. \cite{JP})
$$\gamma_D\le B_D.$$ On the other hand, there
exists a constant $c_n>0,$ depending only on $n$ such that for any
convex domain  $D\subset\Bbb C^n,$ containing no complex line, the
following inequality holds (see \cite{NP2}):
$$B_D\le c_n\gamma_D.$$

This fact extends to any $\Bbb C$-convex domain as the following
theorem shows.

\begin{theorem}\label{8} There exists a constant
$c_n>0,$ depending only on $n,$ such that for any $\Bbb C$-convex
domain $D\subset\Bbb C^n$, containing no complex lines,\footnote
{Under the given assumptions $D$ is biholomorphic to a bounded
domain (cf. \cite {NPZ}) and hence $B_D$ is well-defined.} one has
that
$$B_D(z;X)d_D(z,X)\le c_n,\quad z\in D, X\in(\CC^n)_\ast.$$
In particular, by Proposition \ref{1},
$$\frac{\kappa_D}{4}\le B_D\le 4c_n\gamma_D.$$
\end{theorem}

To prove Theorem \ref{8}, we shall need a lower geometrical
estimates for the Bergman kernel. For this, similarly to the
convex case (see \cite{NP2}; see also \cite{Hef1,Hef2,Con} and
compare with \cite{Che,McN1,McN2}), we introduce the following
geometrical objects related to an arbitrary domain $D\subset\Bbb
C^n,$ containing no complex lines.

For $z^0\in D=:D_0\subset\CC^n=:H_0$ define
$d_{1,D}(z^0):=\dist(z^0,\partial D)=d_D(z^0)$. Fix an
$a^1\in\partial D$ such that $||a^1-z^0||=d_{1,D}(z^0)$. Let
$l_1=z^0+V_1$ be the complex line passing through $z^0$ and $a^1$.
Let $H_1:=V_1^\perp$ be the $(n-1)$-dimensional complex space
orthogonal to $V_1$. Set $D_1:=D_0\cap(z^0+H_1)$ and
$d_{2,D}(z^0):=\dist_{z^0+H_1}(z^0,\partial_{z^0+H_1} D_1)$. Then
fix a point $a^2\in\partial_{z^0+H_1}(D_1)$ with
$\|a^2-z^0\|=d_{2,D}(z^0)$. Denote by $l_2=z^0+V_2$ the complex
line through $z^0$ and $a^2$. Note that $V_2\subset V_1^\perp$.
Put $H_2:=V_2^\perp\cap H_1$ and define $D_2:=D_1\cap(z^0+H_2)$.
Continuing the previous procedure we are led to an orthonormal
basis (arising from the complex lines $l_1,\dots,l_n$), positive
numbers $d_{1,D}(z^0),\dots,d_{n,D}(z^0)$ and points
$a^1,\dots,a^n$ with $a^j\in\partial_{z^0+H_{j-1}}D_{j-1}$ and
$\|a^j-z^0\|=d_{j,D}(z^0)$.

Set $$p_D(z^0):= d_{1,D}(z^0)\cdots d_{n,D}(z^0).$$

Using these numbers we get the following estimates for the Bergman
kernel.

\begin{theorem}\label{9} Let $D\subset\Bbb C^n$ be a $\Bbb C$-convex
domain containing no complex lines. Then
$$\frac{1}{(16\pi)^n}\le
K_D(z)p_D^{2}(z)\le\frac{(2n)!}{(2\pi)^n}.$$
\end{theorem}

Recall that the constant 16 can be replaced by 4 in the case of
convex domains (see \cite{NP2}).

The next result extends earlier ones treating convex domains of
finite type (cf. \cite{Che,McN2}) and the proof here is easier and
pure geometrical. Take a vector $X\in\Bbb C^n.$ For any point
$z\in D$, decompose $X$ w.r.t to the orthogonal basis mentioned
above, i.e. $X=(X_1(z),\dots,X_n(z)).$

Then the following result is a consequence of Proposition \ref{1}
and Theorem \ref{8} .

\begin{proposition}\label{10} There exists a constant
$c_n>1,$ depending only on $n,$ such that for any $\Bbb C$-convex
domain $D\subset\Bbb C^n$, containing no complex lines, one has
that
$$c_n^{-1}\le F_D(z;X)\left(\sum_j^n\frac{|X_j(z)|}{d_{j,D}(z)}\right)^{-1}\le c_n,$$
where $F_D$ denotes any of the Carath\'eodory, Kobayashi or
Bergman metrics.
\end{proposition}

This result is in the spirit of Corollary \ref{14}.
\smallskip

{\bf Remark.} Proposition \ref{3} and Corollary \ref{5}  hold for
the Bergman metric, if the domain contains no complex lines. (In
fact, then Proposition \ref{3} transports the main result in
\cite{Her} and a result in \cite{NP1} from the convex to the $\Bbb
C$-convex case). Moreover, these and the other results for the
Bergman kernel and metric have local versions on bounded
pseudoconvex domains due to the localization principle for the
Bergman invariants (cf. \cite{JP}).

\section{Proofs}

\begin{proof}[Proof of Proposition \ref{1}]
The upper bound is trivial and holds for any domain $D$, since it
contains the disc with center $z$ and radius $d_D(z,X)$ in
direction $X.$

To prove the lower bound, we may assume that $||X||=1.$ Denote by
$l$ the complex line trough $z$ in direction $X$ and choose $a\in
l\cap\partial D$ such that $||z-a||=d_D(z,X).$ Consider a complex
hyperplane $H$ through $a$ such that $D\cap H=\emptyset$ and
denote by $G$ the projection of $D$ onto $l$ in direction $H.$
Note that $G$ is a simply connected domain (cf. \cite{APS,Hor}),
$a\in\partial G$ and $d_G(z)=||z-a||.$ It remains to apply the
Koebe theorem to get that
$$\gamma_D(z;X)\ge\gamma_G(z;1)\ge\frac{1}{4 d_G(z)}.$$
\end{proof}

Many of the next proofs will be based on the following geometric
property of weakly linearly convex domains (see also \cite{ZZ} and
(for the finite type case) \cite{Con}).

\begin{lemma}\label{11} Assume that a weakly linearly
convex domain $G\subset\Bbb C^n$ contains the unit disc $\Bbb D_j$
in the $j$-th complex coordinate line for any $j=1,\dots,n$. Then
$G$ contains the convex hull of $\bigcup_{j=1}^n \Bbb D_j,$ i.e.
$$E:=\{z\in\Bbb C^n:\sum_{j=1}^n|z_j|<1\}\subset G.$$
\end{lemma}

\begin{proof} For any $\eps\in (0,1)$
there is $\delta>0$ such that
$$
X_\eps:=\bigcup_{j=1}^n\Big(\delta\Bbb
D\times\dots\times\delta\Bbb D\times\eps\Bbb D\times \delta\Bbb
D\times\dots\times\delta\Bbb D\Big)\subset G.$$ Recall that
$$
\widehat X_\eps\subset G,
$$
where $\widehat X_\eps$ is the smallest linearly convex set
containing $X_\eps$. Moreover,
$$
\widehat X_\eps=\{z\in\CC^n:\forall b\in\CC^n:<z,b>=1\ \exists
a\in X_\eps:<a,b>=1\}.
$$
(cf. \cite{APS,Hor}). Then $\widehat X_\eps$ is a balanced domain
and, therefore, convex (see \cite{NPZ}). Hence,
$$E_\eps:=\{z\in\Bbb C^n:\sum_{j=1}^n|z_j|<\eps\}\subset\widehat X_\eps\subset
G,\quad \eps\in (0,1),$$ which proves Lemma \ref{11}.
\end{proof}

\noindent{\bf Remark.} The same argument implies that $G$ contains
the convex hull of any balanced domain lying in $G.$ In
particular, the maximal balanced domain lying in $G$ is convex
(see also \cite{ZZ}).

\begin{proof}[Proof of Proposition \ref{3}] (i) Assuming the contrary,
we may find an $r>0$ and sequences $D\supset(z_j)_j,\; z_j\to a,\
\CC^n\supset(X_j)_j,\;X_j\to X\not\in L_a$ such that $\ds
\gamma_D(z_j;X_j)\le\frac{1}{4r}.$ Note that, by Proposition 1,
$d_D(z_j;X_j)\ge r$ (this is trivial if $D$ contains the complex
line through $z_j$ in direction $X_j$). Then
$\Delta_{X_j}(a,r)\subset D_r=D\cap\Bbb B_n(a,2r)$ for any large
$j.$ Note that $D_r$ is a (weakly) linearly convex open set. It is
easy to see that $D_r$ is taut, i.e. the family $\mathcal O(\Bbb
D,D_r)$ is normal (cf. \cite{NPTZ}). Hence
$\Delta_X(a,r)\subset\partial D;$ a contradiction.

(ii) Recall that $\partial D$ is $C^1$-smooth. Therefore, for any
two linearly independent vectors $X,Y\in L_a,$ we may find a
neighborhood $U$ of $a$ and a number $\varepsilon>0$ such that
$\Delta_X(z,\varepsilon)\subset D$ and
$\Delta_Y(z,\varepsilon)\subset D$ for $z\in D\cap U\cap\Lambda.$
It follows by Lemma \ref{11} that
$\Delta_{X+Y}(z,\varepsilon')\subset D$ for some $\varepsilon'>0.$
We get as in (i) that $\Delta_{X+Y}(a,\varepsilon')\subset\partial
D.$ Therefore, $L_a$ is a linear space.

Then, choosing a basis in $L_a$ and applying Lemma \ref{11}, we
see that there are a neighborhood $U$ of $a$ and a number $c>0$
such that $\Delta_X(z,c)\subset D$ for any $z\in D\cap
U\cap\Lambda$ and any unit vector $X\in L_a.$ Now the desired
estimates follow by Proposition \ref{1}.
\end{proof}

\begin{proof}[Proof of Proposition \ref{4}] We may assume that $\text{Re}(z_1)<0$
is the inner normal direction to $\partial D$ at $a=0.$ Let
$r(z)=\text{Re}(z_1)+o(|z_1|)+\rho('z)$ be a smooth defining
function of $D$ near $0.$

For any small $\delta>0$ we have that $\delta=d_D(\delta_n),$
where $\delta_n=(-\delta,'0).$ Set
$L_\delta(\zeta)=-\delta_n+\zeta X,$ $\zeta\in\CC^n.$

We shall consider two cases.

1. $l_{a,X}=1.$ This means that $X_1\neq 0.$ Then
$r(L_\delta(\zeta))=-\delta+\text{Re}(\zeta X_1)+o(|\zeta|).$ It
follows that $L_\delta(\zeta)\in D$ if
$|\zeta|<\frac{\delta}{2|X_1|}$ and $\delta$ is small enough. This
proves the left-hand side inequality.

The opposite inequality follows by the inequality
$r(L_\delta(2\delta/X_1))>0$ which holds for any small $\delta>0.$

2. $l_{a,X}\ge 2.$ This means that $X_1=0.$ Then
$r(L_\delta(\zeta))=-\delta+\rho(\zeta 'X ).$ Since $\rho(\zeta
'X)\le c|\zeta|^l$ for some $c>0,$ we conclude that
$L_\delta(\zeta)\in D$ if $c|\zeta|^l<\delta.$ This implies the
left-hand side inequality.

To prove the opposite inequality, we have to find $c_1>0$ such
that for any small $\delta>0$ there is $\zeta$ with
$|\zeta|^l=c_1^{-1}\delta$ and $\rho(\zeta'X)\ge\delta.$ Since $D$
is (weakly) linearly convex, it follows that
$\rho(\zeta'X)=h(\zeta)+o(|\zeta|^l)\ge 0,$  where
$$h(\zeta)=\sum_{j+k=l}a_{jk}\zeta^j\overline{\zeta}^k\not\equiv 0.$$
Then the homogeneity of $h$ implies that $h\ge 0.$ Moreover, since
$h\not\equiv 0$ we may find a $\zeta$ with $|\zeta|=1$ and
$h(\zeta)>c_1$ for some $c_1>0.$ Then the constant $c_1$ does the
job for any small $\delta>0.$
\end{proof}

\begin{proof}[Proof of Proposition \ref{6}]
The inequality $l_a\le m_a$ is trivial. To prove the opposite one,
we may assume that $l_a<\infty.$ It follows from Propositions
\ref{1} and \ref{4} that
$$\liminf_{D\cap n_a\ni z\to a}\gamma_D(z;X)d^{1/l_a}\ge c_X>0.$$
Hence, $m_a\le l_a$ by Corollary 2 in \cite{Yu4} (in fact,
$\limsup$ instead of $\liminf$ above is sufficient).
\end{proof}

\begin{proof}[Proof of Proposition \ref{7}]
We shall use the same notations as in the proof of Proposition
\ref{4}. It is enough to show that if $\varphi:\Bbb D\to\partial
D$ is a non-trivial analytic disc with $\varphi(0)=0,$ then
$L_a\neq\{0\}.$ Since $\partial D$ is smooth near $a,$ it follows
that there is a $c>0$ such that
$\varphi_\delta(\zeta)=-\delta_n+\varphi(\zeta)\in D$ if
$\delta<c$ and $|\zeta|<c.$ Let $m=\ord_0\varphi$ and
$X=\frac{\varphi^{(m)}(0)}{m!}.$ Denoting by $\kappa^{(m)}_D$ the
Kobayashi metric of order $m$ (cf. \cite{Yu4} for this notion), it
follows that $\kappa^{(m)}_D(\delta_n;X)\le1/c.$ Since
$\gamma_D\le \kappa^{(m)}_D,$ we get as in the proof of
Proposition 3 (i) that $\Delta_X(a,c/4)\subset\partial D.$
\end{proof}

\begin{proof}[Proof of Proposition \ref{12}] The proof can be done following
line by line the proofs in \cite{Yu1}. We only point out how the
replace the arguments there that use convexity. We may assume that
$D$ is a $\Bbb C$-convex domain and $a=0.$ Following the notation
from Proposition \ref{4}, let $r(z)=Re(z_1)+o(|z_1|)+\rho('z)$ be
a defining function of $D$ which is smooth near $0.$

page 841: Let $X,Y\subset\CC^{n-1}$ be such that $\rho(\zeta X)\le
C|\zeta|^m$ and  $\rho(\zeta Y)\le C|\zeta|^m.$ We have to show
that $\rho(\zeta(X+Y)/2)\le C|\zeta|^m.$ For this, fix $\zeta\neq
0$ and take $\delta=C|\zeta^m|.$ Then
$\Delta_X(\delta_n,|\zeta|)\subset D,$
$\Delta_Y(\delta_n,|\zeta|)\subset D$ and hence
$\Delta_{(X+Y)/2}(\delta_n,|\zeta|)\subset D$ by Lemma \ref{11}.
This implies the desired inequality.

We may do the same to get the formula (2.13) on page 845.

Our Proposition \ref{6} is an extension of Theorem C which is
invoked on page 845.

It remains to show Proposition 2 on page 843. Let $k_2,\dots,k_n$
be even integers such that $\rho(\zeta e_j)\le C|\zeta|^{k_j}$ for
any $j=2,\dots,n.$ It is enough to prove that $D^\alpha{\overline
D}^\beta\rho(0)=0$ for any $n$-tuples
$\alpha=(\alpha_2,\dots,\alpha_n)$ and
$\beta=(\beta_2,\dots,\beta_n)$ of non-negative integers with
$w_{\alpha,\beta}=\sum_{j=2}^n\frac{\alpha_j+\beta_j}{k_j}<1.$
Since $\Delta_{e_j}(C\delta_n,\delta^{1/k_j})\subset D$ for any
$\delta>0,$ it follows by Lemma 11 that $\rho('z/n)<C\delta$ for
any $z$ with $|z_j|^{k_j}<\delta.$ In particular, if
$\rho_t('z)=\rho(t^{1/k_2}z_2,\dots,t^{1/k_n}z_n),$ $t>0,$ then
\begin{equation}\label{17}
0\le \rho_t('z/n)<Ct,\quad 'z\in\Bbb D^{n-1}.
\end{equation}
Let now $s=\min\{w(\alpha,\beta):D^\alpha{\overline
D}^\beta\rho(0)\neq0\}.$ Then
$$\lim_{t\to 0}t^{-s}\rho_t('z)=\sum_{w(\alpha,\beta)=s}D^\alpha{\overline
D}^\beta\rho(0)'z^\alpha\overline{'z}^\beta$$ locally uniformly in
$'z.$ Assuming $s<1,$ the inequality \ref{17} implies that the
last polynomial vanishes, a contradiction.
\end{proof}

\begin{proof}[Proof of Proposition \ref{18}] Let first
$L_a\neq\{0\}.$ This means that $\Delta_X(a,r)\subset\partial D$
for some $r>0$ and $X\in(\CC^n)_\ast.$ By convexity,
$\Delta_X(c,r/2)\subset D$ for any $c=ta+(1-t)b$ if $b\in D$ and
$t\in(0,1/2].$ Now the maximum principle implies that $a$ is not a
peak point.

Let now $L_a=\{0\}.$ We may assume that $a=0$ and
$D\subset\{z\in\CC^n:\text{Re}(z_1)<0\}.$ Then $e^{z_1}$ is an
entire weak peak function for $\overline{D}$ at $0.$ Setting
$H=\{z\in\CC^n:\text{Re}(z_1)<0\}.$ It follows that implies
$\text{supp}\mu\subset D_1=\partial D\cap H$ for any representing
measure $\mu$ for $0$ w.r.t. $A(D).$ Since $L_0=\{0\},$ it follows
that $0$ is a boundary point of the convex set $D_1.$ Then there
exists an entire function which is a weak peak function for $D_1$
at $0$ (we need such a function function to be in $A(D)$). We get
as above that $\text{supp}\mu$ is contained in some $(n-2)$
dimensional space. Repeating this procedure, it follows that
$\text{supp}\mu\subset\partial D\cap l,$ where $l$ is a complex
line. Since $0$ is a boundary point of the last convex set, then
there is an entire function which is a peak function for $\partial
D\cap l$ at $0.$ So $\text{supp}\mu=\{0\},$ i.e. 0 is a peak point
w.r.t. $A(D)$ (cf. \cite{Gam}).
\end{proof}

\begin{proof}[Proof of Theorem \ref{9}.] We first prove the lower
bound. Fix $z^0\in D$. Using a translation and then successive
rotations we may assume (see the description of the numbers
$d_{j,D}$) that $z^0=0$, $H_j=\{0\}\times \CC^{n-j}$,
$j=1,\dots,n-1$, and
$a^j=(0,a^j_j,0)\in\CC^{j-1}\times\CC\times\CC^{n-j}$ with
$d_{j,D}(z^0)=|a_j^j|$.

Recall that $D$ is $\CC$-convex. Therefore, there exist affine
hyperplanes $a^j+W_j$ through $a^j$ which do not intersect $D$.
Note that $W_1\cap H_1$ is orthogonal to $a^2$, i.e. $W_1\cap
H_1\subset\{0\}\times\CC^{n-2}$. Hence $W_1$ is given by the
equation $\alpha_{2,1}z_1+z_2=0$. Moreover, using a similar
argument, the equations for $W_j$, $j=1,\dots,n-1$, are the
following ones:
$$
\alpha_{j,1}z_1+\cdots+\alpha_{j,j}z_j+z_{j+1}=0.
$$
Let $F:\CC^n\to\CC^n$ be the linear mapping given by the matrix
$A$ whose rows are given by the vectors
$(\alpha_{j,1},\dots,\alpha_{j,j},1,0,\dots,0),$ $j=0,\dots,n-1.$
Define $G=F(D)$ and observe that $G$ is again $\CC$-convex. Note
that $K_D(0)=K_G(0)$ since $\det A=1.$ Finally, put
$G_j:=\pi_j(G)$, where $\pi_j$ is the projection onto the $j$-th
coordinate axis. Then (see \cite{APS}) $G_j$ is a simply connected
domain, $j=1,\dots,n,$ and $G\subset G_1\times\cdots\times G_n$.
Hence
$$
K_D(0)\geq K_{G_1\times\cdots\times G_n}(0)=K_{G_1}(0)\cdots
K_{G_n}(0).
$$
Since $G_j$ is simply connected, using the Koebe theorem we get
$$
\sqrt{\pi K_{G_j}(0)}=\gamma_{G_j}(0;1)\geq \frac{1}{4d_{G_j}(0)}.
$$
Note that $F(a^j)\in\partial G$, its $j$-th coordinate is $a_j^j,$
and the affine hyperplane $\{z\in\Bbb C^n:z_j=a_j^j\}$ does not
intersect $G$. Hence $a_j^j\in\partial G_j$; in particular,
$d_{j,D}(z^0)=|a_j^j|\geq d_{G_j}(0)$, which finally gives the
lower bound.

To show the upper bound, consider the dilatation of coordinates
$$\Phi(z)=(z_1/d_{1,D}(z^0),\dots,z_n/d_{n,D}(z^0))$$ and
set $\tilde G=\Phi(D).$ Hence
$$K_D(z^0)=\frac{K_{\tilde G}(0)}{p_D^2(z^0)}.$$

Then the upper bound follows from Lemma \ref{11} and the following
formula (cf. \cite{JP,NP1}):
$$K_E(0)=\frac{(2n)!}{(2\pi)^n}.$$
\end{proof}

\begin{proof}[Proof of Theorem \ref{8}.] The proof can be done following
line by line the proof of Theorem 2 in \cite{NP2} and using
Theorem \ref{9} and Lemma \ref{11}. For convenience of the reader,
we include a complete proof.

We shall use the geometric constellation in the proof of Theorem
\ref{9}. Let $X\in(\CC^n)_\ast$ and fix $k\in J:=\{j:X_j\neq 0\}$.
Then
$$\Psi_k(z):=(z_1-\frac{X_1}{X_k}z_k,\ldots,z_{k-1}-\frac{X_{k-1}}{X_k}z_{k},
z_k,z_{k+1}-\frac{X_{k+1}}{X_k}z_k,\ldots,z_n-\frac{X_n}{X_k}z_k)$$
is a linear mapping with jacobian equal to 1 and $Y^k:=\Psi_k(X)=
(0,\ldots,0,X_k,0,\ldots,0)$. Let $\Delta_j$ be the disc in the
$j$-th coordinate plane with center at $0$ and radius $d_{j,D}(0)$
if $j\neq k$, and $d'_k:=|X_k|d_D(0,X)$ if $j=k$. Then
$\Delta_j\subset D_k:=\Psi_k(D)$ and, by Lemma \ref{11},
$$D_k\supset E_k=\{z\in\CC^n:\frac{|z_k|}{d'_k}+\sum_{j=1,j\neq
k}^n\frac{|z_j|}{d_j}<1\}.$$ Hence $$M_D(0;X)=M_{D_k}(0;Y^k)\le
M_{E_k}(0;Y^k)=C\frac{d_{k,D}(0)}{|X_k|p_D(0)d_D^2(0,X)},$$ where
$C_n:=M_E(0;e_1)=\sqrt{\frac{(2(n+1))!}{6(2\pi)^n}}$ (cf.
\cite{NP2}) and $e_1$ is the first basis vector. Applying the
lower bound in Theorem \ref{9},  we obtain that
\begin{equation}\label{15}B_D(0;X)=\frac{M_D(0;X)}{\sqrt{K_D(0)}}\le
\frac{c'_nd_{k,D}(0)}{|X_k|d_D^2(0,X)}, \quad 1\le k\le n,
\end{equation}
where
$c'_n=(4\sqrt{\pi})^nC_n=2^n\sqrt{\frac{2^{n-1}(2(n+1))!}{3}}.$ It
remains to apply Lemma \ref{11} to get that
\begin{equation}\label{16}
\frac{1}{d_D(0,X)}\le\sum_j^n\frac{|X_j(z)|}{d_{j,D}(z)}
\end{equation}
and then to choose $c_n=nc'_n.$
\end{proof}

\begin{proof}[Proof of Proposition \ref{10}.] It follows by
(\ref{15}) and the inequality
$$B_D(z;X)\ge\frac{1}{4d_D(z,X)}$$
that
$$\frac{|X_j(z)|}{d_{j,D}(z)}\le\frac{4c_n'}{d_D(z)}.$$
Hence,
\begin{equation}\label{19}
\frac{1}{d_D(z,X)}\le\sum_j^n\frac{|X_j(z)|}{d_{j,D}(z)}\le\frac{4c_n}{d_D(z,X)},
\end{equation}
where $c_n=nc_n'.$ Then (\ref{15}) and (\ref{16}) imply that
$$(16c_n)^{-1}\le F_D(z;X)\left(\sum_j^n\frac{|X_j(z)|}{d_{j,D}(z)}\right)^{-1}\le
c_n.$$
\end{proof}

\noindent{\bf Remarks.} (a) In \cite{Che,McN2}, the numbers
$d_{j,D}(z)$ are replaced by other numbers $\tilde d_{j,D}(z)$ in
the finite type convex case. Note that $\tilde
d_{1,D}(z)=d_{1,D}(z)$ and $\tilde d_{j,D}(z)\ge \tilde
d_{j+1,D}(z),$ $2\le j\le n-1$ (in contrast to $d_{j,D}(z)\le
d_{j+1,D}(z)$). The inductive definition of $\tilde d_{j,D}(z)$ is
similar to that of $d_{j,D}(z)$ but in any step $j\ge 2$ the
number $d_{j,D}(z)$ is the radius of the largest (not the
smallest!) disc in the respective $(n-j+1)$-dimensional set.
However, one can show the the respective supporting hyperplanes
$\tilde W_j$ have the same equations as $W_j$ when $z$ is near
$\partial D.$ Then the above approach allows us to get the same
estimates as in Theorem \ref{9} and Proposition \ref{10} in terms
of these numbers and the respective coordinates. In particular, it
leads to (\ref{19}) also for the $\tilde d_{j,D}$-situation.

(b) Assume that a domain $D\subset\CC^n$ is smooth and weakly
linearly convex near a boundary point $a$ of finite type $m.$ Then
$\tilde d_{j,D}(z)\le (d_D(z))^{1/m}$ for $j=1,\dots,n$ (see
\cite{Blu}). Since $a$ is a local peak point for $D$ at $a,$ it
follows that there is a neighborhood $U$ of $a$ and a constant
$c>0$ such that
$$\kappa_D(z;X)\ge\frac{c||X||}{(d_D(z))^{1/m}},\quad z\in D\cap U;$$
the same estimate holds for $B_D$ if $D$ is pseudoconvex (not
necessary bounded - use e.g. localization results in \cite{Nik2}).

(c) Finally, note that there is a number $c_n>1$ depending only on
$n$ such that
$$c_n^{-1}\le \frac{d_{j,D}(z)}{\tilde d_{n-j,D}(z)}\le c_n,\quad z\in D,j=1,\dots,n,$$
for any $z$ near the boundary of any smooth convex domain
$D\subset\CC^n$ containing no complex lines \footnote{It will be
interesting to have a pure geometric proof of this inequality and
for (\ref{19}) with best possible constants.} (see also
\cite{Hef2}).

In fact, a more general statement is true. Let $\{p_1,\dots,p_n\}$
and $\{q_1,\dots,q_n\}$ be orthonormal bases in $\CC^n.$ Let
$a_1,\dots,a_n$ and $b_1\dots,b_n$ be increasing sequences of
positive numbers. Assume that there is $c>1$ such that
$$c^{-1}\le\frac{\sum_{j=1}^na_j|\langle
X,p_j\rangle|}{\sum_{j=1}^nb_j|\langle X,q_j\rangle|}\le
c\quad\text{ for any }X\in(\CC^n)_\ast.$$ Then
$$c'^{-1}\le\frac{a_j}{b_j}\le c',\quad j=1,\dots,n,$$ where
$c'=n!c.$

For this, observe that expanding the determinant of the matrix of
the unitary transformation between the bases, it follows that
$$\prod_{j=1}^n|\langle p_j,q_{\sigma(j)}\rangle|\ge\frac{1}{n!}$$
for some permutation $\sigma$ of $\{1,\dots,n\}.$ In particular,
$|\langle p_j,q_{\sigma(j)}\rangle|\ge 1/n!.$ Then the given
condition implies that
$$c'^{-1}\le\frac{a_j}{b_{\sigma(j)}}\le c'.$$
Assume now that $a_k>c'b_k$ for some $k.$ Using the monotonicity
and the inequality $a_j\le c' b_{\sigma(j)},$ it follows that
$\sigma(j)>k$ for any $j\ge k.$ Since $\sigma$ is a permutation of
$\{1,\dots,n\},$ we get a contradiction.

\end{document}